\documentclass[11pt]{article}
\usepackage{amsmath}
\usepackage{amsfonts}
\usepackage{graphicx}

\input{tcilatex}
\begin{document}

\begin{center}
{\Large Convergence of the Bernstein-Durrmeyer operators in variation
seminorm}\bigskip 

\vskip0.2in

\"{O}zlem \"{O}ks\"{u}zer$^{a}$,\ Harun Karsli$^{b}$\ and Fatma Ta\c{s}delen
Ye\c{s}ildal$^{a}$\vskip0.3in

$^{a}$Department of Mathematics, Faculty of Science,

University of Ankara, Ankara, TURKEY

oksuzer@ankara.edu.tr; tasdelen@science.ankara.edu.tr\vskip0.3in

$^{b}$Abant Izzet Baysal University, Faculty of Science and Arts,

Department of Mathematics, Bolu / Turkey

karsli\_h@ibu.edu.tr

\vskip0.3in
\end{center}

{\footnotesize Abstract. The aim of this paper is to study variation
detracting property and convergence in variation of the Bernstein-Durrmeyer
modifications of the classical Bernstein operators in the space of functions
of bounded variation. These problems are studied with respect to the
variation seminorm. Moreover we also study the problem of the rate of
approximation.}\medskip

\medskip

{\footnotesize AMS Subject Classification 2000: 41A25, 41A35\medskip , 41A36}

\begin{center}
\textbf{1. Introduction\medskip }
\end{center}

The classical Bernstein operators, defined for functions $f\in C[0,1]$, are
of the form%
\begin{equation}
\left( B_{n}f\right) \left( x\right) =\sum\limits_{k=0}^{n}f\left( \frac{k}{n%
}\right) p_{n,k}\left( x\right) ~,n\in 
%TCIMACRO{\U{2115} }%
%BeginExpansion
\mathbb{N}
%EndExpansion
\label{be}
\end{equation}%
where $p_{n,k}\left( x\right) :=\binom{n}{k}x^{k}\left( 1-x\right) ^{n-k}$
is the Bernstein basis function, in which the function $f$ can be completely
reconstructed for all $x\in \lbrack 0,1]$ from its sampled value $f(k/n)$
taken at the nodes $k/n$ $(k=0,1,...,n).$\medskip

Let $L_{1}=L_{1}\left[ 0,1\right] $ be the space of all Lebesgue measurable
and integrable complex valued functions defined on $\left[ 0,1\right] $,
endowed with the usual norm. In order to find a positive answer to the
approximation problem for Lebesgue integrable functions defined on the
interval $\left[ 0,1\right] ,$ J. L. Durrmeyer \cite{1} and, independently
Lupas \cite{2} introduced the integral modifications of the well-known
Bernstein polynomials, called Bernstein-Durrmeyer operators. These are
defined as%
\begin{equation}
\left( D_{n}f\right) \left( x\right) =\left( n+1\right)
\sum\limits_{k=0}^{n}p_{n,k}\left( x\right) \int\limits_{0}^{1}p_{n,k}\left(
t\right) f\left( t\right) dt,\text{\ }x\in \left[ 0,1\right] .  \label{1}
\end{equation}%
Here, the sampled value $f(k/n)$ in (\ref{be}) is replaced by%
\begin{equation*}
\left( n+1\right) \int\limits_{0}^{1}p_{n,k}\left( t\right) f\left( t\right)
dt.
\end{equation*}%
In 1981 Derriennic \cite{3} first studied these operators in details. After
that, Bernstein-Durrmeyer operators have been extensively studied by Gonska
and Zhou \cite{4}, Ditzian and Ivanov \cite{5}, and several other
authors.\medskip

The present work is strongly motivated by the paper \cite{6}, in which the
authors have introduced, developed in details and studied the deep
interconnections between.variation detracting property and the convergence
in variation for Bernstein-type polynomials and singular convolution
integrals. After this fundamental study, the convergence in variation
seminorm has become a new research field in the theory of approximation. One
should also note the results obtained in this direction in [\cite{8}-\cite%
{12}].

The aim of the paper is to establish variation detracting property and
convergence in variation of Bernstein-Durrmeyer operators in the space of
functions of bounded variation. The rate of approximation is given with
respect to the variation seminorm.\medskip

\begin{center}
\textbf{2. Notation and Preliminaries}\bigskip
\end{center}

For the notation; let $I$\ be a bounded or unbounded interval. Throughout
the work, $V_{\left[ I\right] }\left[ f\right] $ stands for the total Jordan
variation of the real-valued function $f$\ defined on $I$. We dealt with the
classes $TV\left( I\right) $ and $BV\left( I\right) $ of all the functions
of bounded variation on $I\subset 
%TCIMACRO{\U{211d} }%
%BeginExpansion
\mathbb{R}
%EndExpansion
$, endowed with the seminorm and norm, respectively%
\begin{equation*}
\left\Vert f\right\Vert _{TV\left( I\right) }:=V_{\left[ I\right] }\left[ f%
\right] ,
\end{equation*}%
and%
\begin{equation*}
\left\Vert f\right\Vert _{BV\left( I\right) }:=V_{\left[ I\right] }\left[ f%
\right] +\left\vert f\left( c\right) \right\vert ,
\end{equation*}%
where $c\ $is any fixed point of $I.$ Some interesting properties of the
space $TV(I)$ are presented in \cite{6}.\medskip

In order to obtain a convergence result in the variation seminorm, it is
necessary and important to state the variation detracting property. It was
G. G. Lorentz who showed in his candidate thesis \cite{7} that the operators 
$B_{n}$ have the property%
\begin{equation*}
V_{\left[ 0,1\right] }\left[ B_{n}f\right] \leq V_{\left[ 0,1\right] }\left[
f\right] \text{ \ \ }\left( n\in 
%TCIMACRO{\U{2115} }%
%BeginExpansion
\mathbb{N}
%EndExpansion
\right) ,
\end{equation*}%
called the variation detracting or variation diminishing property, i.e.,
positive linear operators from the space of functions of bounded variation
into itself do not increase the total variation of functions.\medskip

Set $AC\left( I\right) ,$ the space of all absolutely continuous real-valued
functions on $I$, is a closed subspace of $TV\left( I\right) $ with respect
to the convergence induced by the seminorm $\left\Vert f\right\Vert
_{TV\left( I\right) }$. In addition, it is well known that is $\underset{%
n\rightarrow \infty }{\lim }V_{I}\left[ g_{n}-g\right] =0$ for a sequence $%
\left( g_{n}\right) _{n\geq 1}$ in $AC\left( I\right) $, then also $g\in AC%
\left[ 0,1\right] $ and%
\begin{equation*}
V_{I}\left[ g_{n}-g\right] =\int\limits_{I}\left\vert g_{n}^{\prime }\left(
t\right) -g^{\prime }\left( t\right) \right\vert dt.
\end{equation*}%
So, convergence in variation of $\left( g_{n}\right) _{n\geq 1}\subset
AC\left( I\right) $ to $g,$ exactly means the convergence of the derivatives 
$g_{n}^{\prime }$ to $g^{\prime }$\ in the $L_{1}\left( I\right) -$%
norm.\medskip

For the proof of theorems in the following sections, since $\left(
d/dx\right) p_{k,n}\left( x\right) =\left( k-nx\right) p_{k,n}\left(
x\right) /X,$ we calculate the two fundamental representations for
derivative of $\left( D_{n}f\right) \left( x\right) $ as,%
\begin{equation}
\left( D_{n}f\right) ^{\prime }\left( x\right) =\frac{\left( n+1\right) }{X}%
\sum\limits_{k=0}^{n}\left( k-nx\right) p_{n,k}\left( x\right)
\int\limits_{0}^{1}f\left( t\right) p_{n,k}\left( t\right) dt  \label{a}
\end{equation}%
and%
\begin{equation}
\left( D_{n}f\right) ^{\prime }\left( x\right)
=n\sum\limits_{k=0}^{n-1}p_{n-1,k}\left( x\right) \left( n+1\right)
\int\limits_{0}^{1}f\left( t\right) \left[ p_{n,k+1}\left( t\right)
-p_{n,k}\left( t\right) \right] dt  \label{b}
\end{equation}%
where $X=x\left( 1-x\right) .$ In the same proofs we need sum moments for
the operators (\ref{1}). Let us define the sum moments for $r=0,1,2,3,4$,%
\begin{equation}
T_{r,n}\left( x\right) =\sum\limits_{k=0}^{n}k^{r}p_{n,k}\left( x\right) .
\label{2}
\end{equation}%
Then there hold%
\begin{equation}
T_{r,n}\left( x\right) =\left\{ 
\begin{tabular}{ll}
$1$ & $r=0$ \\ 
$nx$ & $r=1$ \\ 
$n\left( n-1\right) x^{2}+nx$ & $r=2$ \\ 
$n\left( n-1\right) \left( n-2\right) x^{3}+3n\left( n-1\right) x^{2}+nx$ & $%
r=3$ \\ 
$%
\begin{array}{c}
n\left( n-1\right) \left( n-2\right) \left( n-3\right) x^{4}+6n\left(
n-1\right) \left( n-2\right) x^{3}+ \\ 
+7n\left( n-1\right) x^{2}+nx%
\end{array}%
$ & $r=4.$%
\end{tabular}%
\right.  \label{mom}
\end{equation}

\begin{center}
\textbf{3. Variation Detracting Property of Bernstein Durrmeyer Operators}$%
\medskip $
\end{center}

In this section, we state the variation detracting property of the
Bernstein-Durrmeyer Operators$\medskip $.

\textbf{Theorem 3.1. }If $f\in TV\left[ 0,1\right] $, then%
\begin{equation}
V_{\left[ 0,1\right] }\left[ D_{n}f\right] \leq V_{\left[ 0,1\right] }\left[
f\right]  \label{14}
\end{equation}%
and 
\begin{equation}
\left\Vert D_{n}f\right\Vert _{BV\left[ 0,1\right] }\leq \left\Vert
f\right\Vert _{BV\left[ 0,1\right] }  \label{15}
\end{equation}%
hold true.\medskip

\textbf{Proof.\ }For convenience we write the Bernstein-Durrmeyer operators
as;%
\begin{equation*}
\left( D_{n}f\right) \left( x\right) =\sum\limits_{k=0}^{n}p_{n,k}\left(
x\right) F_{k,n}
\end{equation*}%
where%
\begin{equation*}
F_{k,n}:=\left( n+1\right) \int\limits_{0}^{1}f\left( t\right) p_{n,k}\left(
t\right) dt.
\end{equation*}%
As in the (\ref{b}), we calculate differentiation of (\ref{1})%
\begin{eqnarray}
\left( D_{n}f\right) ^{\prime }\left( x\right)
&=&\sum\limits_{k=0}^{n}p_{n,k}^{\prime }\left( x\right)
F_{k,n}=\sum\limits_{k=1}^{n}\binom{n}{k}kx^{k-1}\left( 1-x\right)
^{n-k}F_{k,n}  \notag \\
&&-\sum\limits_{k=0}^{n-1}\binom{n}{k}x^{k}\left( n-k\right) \left(
1-x\right) ^{n-k-1}F_{k,n}  \notag \\
&=&n\sum\limits_{k=0}^{n-1}p_{n-1,k}\left( x\right)
F_{k+1,n}-n\sum\limits_{k=0}^{n-1}p_{n-1,k}\left( x\right) F_{k,n}  \notag \\
&=&n\sum\limits_{k=0}^{n-1}p_{n-1,k}\left( x\right) \left[ F_{k+1,n}-F_{k,n}%
\right]  \notag \\
&=&n\sum\limits_{k=0}^{n-1}p_{n-1,k}\left( x\right) \Delta F_{k,n}.
\label{3}
\end{eqnarray}%
Considering the representation (\ref{3}) of\textbf{\ }$\left( D_{n}f\right)
^{\prime }$, one has%
\begin{eqnarray*}
\left\Vert D_{n}f\right\Vert _{TV\left[ 0,1\right] } &=&V_{\left[ 0,1\right]
}\left[ D_{n}f\right] =\int\limits_{0}^{1}\left\vert \left( D_{n}f\right)
^{\prime }\left( x\right) \right\vert dx \\
&\leq &n\sum\limits_{k=0}^{n-1}\left\vert \Delta F_{k,n}\right\vert
\int\limits_{0}^{1}p_{n-1,k}\left( x\right) dx.
\end{eqnarray*}%
Since $n\int\limits_{0}^{1}p_{n-1,k}\left( x\right) dx=1$, we get%
\begin{equation}
\left\Vert D_{n}f\right\Vert _{TV\left[ 0,1\right] }\leq
\sum\limits_{k=0}^{n-1}\left\vert \Delta F_{k,n}\right\vert .  \label{*}
\end{equation}%
Here $\Delta F_{k,n}$ is%
\begin{eqnarray*}
\Delta F_{k,n} &=&\left( n+1\right) \int\limits_{0}^{1}f\left( t\right) 
\left[ p_{n,k+1}\left( t\right) -p_{n,k}\left( t\right) \right] dt \\
&=&\left( n+1\right) \int\limits_{0}^{1}f\left( t\right) \Delta
p_{n,k}\left( t\right) dt.
\end{eqnarray*}%
Since%
\begin{eqnarray*}
p_{n,k}^{\prime }\left( x\right) &=&n\left( p_{n-1,k-1}\left( x\right)
-p_{n-1,k}\left( x\right) \right) \\
&=&-n\left( p_{n-1,k}\left( x\right) -p_{n-1,k-1}\left( x\right) \right) \\
&=&-n\Delta p_{n-1,k-1}\left( x\right) ,
\end{eqnarray*}%
and%
\begin{equation*}
\Delta p_{n-1,k-1}\left( t\right) =\frac{-p_{n,k}^{\prime }\left( t\right) }{%
n}\Longrightarrow \Delta p_{n,k}\left( t\right) =\frac{-p_{n+1,k+1}^{\prime
}\left( t\right) }{n+1},
\end{equation*}%
we get%
\begin{equation}
\left\vert \Delta F_{k,n}\right\vert =\left\vert \left( n+1\right)
\int\limits_{0}^{1}f\left( t\right) \left[ \frac{-p_{n+1,k+1}^{\prime
}\left( t\right) }{n+1}\right] dt\right\vert .  \label{**}
\end{equation}%
So from (\ref{*}) and (\ref{**}), we obtain%
\begin{eqnarray*}
V_{\left[ 0,1\right] }\left[ D_{n}f\right] &=&\sum\limits_{k=0}^{n-1}\left%
\vert \Delta F_{k,n}\right\vert =\sum\limits_{k=0}^{n-1}\left\vert
-\int\limits_{0}^{1}f\left( t\right) p_{n+1,k+1}^{\prime }\left( t\right)
dt\right\vert \\
&=&\sum\limits_{k=0}^{n-1}\left\vert \int\limits_{0}^{1}p_{n+1,k+1}\left(
t\right) f^{\prime }\left( t\right) dt\right\vert \\
&\leq &\sum\limits_{k=0}^{n-1}\int\limits_{0}^{1}p_{n+1,k+1}\left( t\right)
\left\vert f^{\prime }\left( t\right) \right\vert dt \\
&=&\int\limits_{0}^{1}\sum\limits_{k=0}^{n-1}\binom{n+1}{k+1}t^{k+1}\left(
1-t\right) ^{n-k}\left\vert f^{\prime }\left( t\right) \right\vert dt \\
&=&\int\limits_{0}^{1}\sum\limits_{k=1}^{n}\binom{n+1}{k}t^{k}\left(
1-t\right) ^{n+1-k}\left\vert f^{\prime }\left( t\right) \right\vert dt \\
&\leq &\int\limits_{0}^{1}\sum\limits_{k=0}^{n+1}\binom{n+1}{k}t^{k}\left(
1-t\right) ^{n+1-k}\left\vert f^{\prime }\left( t\right) \right\vert dt \\
&=&\int\limits_{0}^{1}\left( t+1-t\right) ^{n+1}\left\vert f^{\prime }\left(
t\right) \right\vert dt=\int\limits_{0}^{1}\left\vert f^{\prime }\left(
t\right) \right\vert dt \\
&=&V_{\left[ 0,1\right] }\left[ f\right] .
\end{eqnarray*}%
The desired estimate (\ref{14}) is now evident.\medskip

Since $\left( D_{n}f\right) \left( 0\right) =\left( n+1\right)
\int\limits_{0}^{1}\left( 1-t\right) ^{n}f\left( t\right) dt$ and $%
\left\Vert f\right\Vert _{BV\left[ I\right] }:=V_{\left[ I\right] }\left[ f%
\right] +\left\vert f\left( 0\right) \right\vert $, relation (\ref{15}) is a
result of (\ref{14}). Indeed%
\begin{eqnarray*}
\left\Vert D_{n}f\right\Vert _{BV\left[ 0,1\right] } &=&V_{\left[ 0,1\right]
}\left[ D_{n}f\right] +\left\vert \left( D_{n}f\right) \left( 0\right)
\right\vert \\
&\leq &V_{\left[ 0,1\right] }\left[ f\right] +\left\vert \left( n+1\right)
\int\limits_{0}^{1}\left( 1-t\right) ^{n}f\left( t\right) dt\right\vert
\end{eqnarray*}%
Since $f\in TV\left[ 0,1\right] $ and $\left\vert \left( n+1\right)
\int\limits_{0}^{1}\left( 1-t\right) ^{n}f\left( t\right) dt\right\vert
=\left\vert f\left( 0\right) \right\vert \leq \left\vert f\left( c\right)
\right\vert $ where $c$ is any fixed point of $\left[ 0,1\right] $, we get%
\begin{equation*}
\left\Vert D_{n}f\right\Vert _{BV\left[ 0,1\right] }\leq \left\Vert
f\right\Vert _{BV\left[ 0,1\right] }.
\end{equation*}%
Thus, the proof of the theorem is complete.\medskip

\begin{center}
\textbf{4. Rate of Approximation in }$TV$\textbf{-norm}\medskip
\end{center}

This section is deal with the rates of approximation $D_{n}g$ to $g$ in the
variation seminorm.\medskip

\textbf{Theorem 4.1.} Let $g^{\prime \prime }\in AC\left[ 0,1\right] $, then%
\begin{equation*}
V_{\left[ 0,1\right] }\left[ D_{n}g-g\right] \leq \frac{2\left( C+1\right) }{%
\sqrt{n}}\left\{ V_{\left[ 0,1\right] }\left[ g\right] +V_{\left[ 0,1\right]
}\left[ g^{\prime \prime }\right] \right\}
\end{equation*}%
where $C>1$ is a constant.\medskip

\textbf{Proof. }By Taylor's formula with integral remainder term, one has%
\begin{equation*}
g\left( t\right) =g\left( x\right) +\left( t-x\right) g^{\prime }\left(
x\right) +\left( t-x\right) ^{2}\frac{g^{\prime \prime }\left( x\right) }{2}+%
\frac{1}{2}\int\limits_{0}^{t-x}\left( t-x-u\right) ^{2}g^{\prime \prime
\prime }\left( x+u\right) du.
\end{equation*}%
Substituting $x+u=v$, it is easily reached that%
\begin{eqnarray}
g\left( t\right) &=&g\left( x\right) +\left( t-x\right) g^{\prime }\left(
x\right) +\left( t-x\right) ^{2}\frac{g^{\prime \prime }\left( x\right) }{2}+%
\frac{1}{2}\int\limits_{x}^{t}\left( t-v\right) ^{2}g^{\prime \prime \prime
}\left( v\right) dv.  \notag \\
&&  \label{4}
\end{eqnarray}%
If we apply the operator $D_{n}^{\prime }$ to both sides of equality (\ref{4}%
), one has from (\ref{a})%
\begin{eqnarray}
\left( D_{n}g\right) ^{\prime }\left( x\right) &=&A_{0,n}\left( x\right)
g\left( x\right) +A_{1,n}\left( x\right) g^{\prime }\left( x\right)
+A_{2,n}\left( x\right) g^{\prime \prime }\left( x\right) +\left(
R_{n}g\right) \left( x\right)  \notag \\
&&  \label{5}
\end{eqnarray}%
where%
\begin{equation*}
A_{j,n}\left( x\right) =\frac{n+1}{Xj!}\sum\limits_{k=0}^{n}\left(
k-nx\right) p_{n,k}\left( x\right) \int\limits_{0}^{1}\left( t-x\right)
^{j}p_{n,k}\left( t\right) dt\text{ \ \ \ }\left( j=0,1,2\right)
\end{equation*}%
and the remainder term is given by%
\begin{equation*}
\left( R_{n}g\right) \left( x\right) =\frac{n+1}{2X}\sum\limits_{k=0}^{n}%
\left( k-nx\right) p_{n,k}\left( x\right) \int\limits_{0}^{1}\left[
\int\limits_{x}^{t}\left( t-v\right) ^{2}g^{\prime \prime \prime }\left(
v\right) dv\right] p_{n,k}\left( t\right) dt
\end{equation*}%
\begin{equation}
\leq \frac{n+1}{2X}\sum\limits_{k=0}^{n}\left( k-nx\right) p_{n,k}\left(
x\right) \int\limits_{0}^{1}\left( t-x\right) ^{2}\left[ \int%
\limits_{x}^{t}g^{\prime \prime \prime }\left( v\right) dv\right]
p_{n,k}\left( t\right) dt.  \label{6}
\end{equation}%
Calculating (\ref{2}), (\ref{mom}) and using the property of the binomial
coefficients, Beta and Gamma functions, we obtain%
\begin{eqnarray*}
A_{0,n}\left( x\right) &=&\frac{1}{X}\left[ T_{1,n}\left( x\right)
-nxT_{0,n}\left( x\right) \right] =0, \\
A_{1,n}\left( x\right) &=&\frac{1}{(n+2)X}\left[ T_{2,n}\left( x\right)
+\left( 1-nx\right) T_{1,n}\left( x\right) -nxT_{0,n}\left( x\right) \right]
=\frac{n}{n+2}
\end{eqnarray*}%
and%
\begin{eqnarray*}
A_{2,n}\left( x\right) &=&\left[ \frac{T_{3,n}\left( x\right) +3(1-\left(
n+2\right) x)T_{2,n}\left( x\right) +\left( 2-\left( 5n+6\right) x+2n\left(
n+3\right) x^{2}\right) T_{1,n}\left( x\right) }{2(n+2)\left( n+3\right) X}%
\right. \\
&&\left. \frac{+\left( 2-\left( 5n+6\right) x+2n\left( n+3\right)
x^{2}\right) T_{1,n}\left( x\right) -2nx\left( 1-\left( n+3\right) x\right)
T_{0,n}\left( x\right) }{2(n+2)\left( n+3\right) X}\right] \\
&=&\frac{2n\left( 1-2x\right) }{\left( n+2\right) \left( n+3\right) }.
\end{eqnarray*}%
So by (\ref{4}), (\ref{5}) and (\ref{6}), we have%
\begin{equation}
\left( D_{n}g\right) ^{\prime }\left( x\right) =\frac{n}{n+2}g^{\prime
}\left( x\right) +\frac{2n\left( 1-2x\right) }{\left( n+2\right) \left(
n+3\right) }g^{\prime \prime }\left( x\right) +\left( R_{n}g\right) \left(
x\right) .  \label{11}
\end{equation}%
In order to estimate the integration domain of the double integral in the
remainder term (\ref{6}), we divide the summation into four different sums
as following ($\left[ x\right] $ denotes the integer part of $x$)%
\begin{equation}
\left( R_{n}g\right) \left( x\right) =\sum\limits_{i=1}^{4}\left(
B_{i,n}g\right) \left( x\right) \equiv \sum\limits_{i=1}^{4}B_{i,n}g
\label{12}
\end{equation}%
where%
\begin{eqnarray*}
B_{1,n}g &=&\frac{1}{2X}\sum\limits_{k=0}^{\left[ nx\right] }\left(
k-nx\right) p_{n,k}\left( x\right) \left( n+1\right)
\int\limits_{0}^{k/n}\left( t-x\right) ^{2}\left[ \int\limits_{x}^{t}g^{%
\prime \prime \prime }\left( v\right) dv\right] p_{n,k}\left( t\right) dt, \\
B_{2,n}g &=&\frac{1}{2X}\sum\limits_{k=0}^{\left[ nx\right] }\left(
k-nx\right) p_{n,k}\left( x\right) \left( n+1\right)
\int\limits_{k/n}^{x}\left( t-x\right) ^{2}\left[ \int\limits_{x}^{t}g^{%
\prime \prime \prime }\left( v\right) dv\right] p_{n,k}\left( t\right) dt, \\
B_{3,n}g &=&\frac{1}{2X}\sum\limits_{k=\left[ nx\right] +1}^{n}\left(
k-nx\right) p_{n,k}\left( x\right) \left( n+1\right)
\int\limits_{x}^{k/n}\left( t-x\right) ^{2}\left[ \int\limits_{x}^{t}g^{%
\prime \prime \prime }\left( v\right) dv\right] p_{n,k}\left( t\right) dt,
\end{eqnarray*}%
and%
\begin{equation*}
B_{4,n}g=\frac{1}{2X}\sum\limits_{k=\left[ nx\right] +1}^{n}\left(
k-nx\right) p_{n,k}\left( x\right) \left( n+1\right)
\int\limits_{k/n}^{1}\left( t-x\right) ^{2}\left[ \int\limits_{x}^{t}g^{%
\prime \prime \prime }\left( v\right) dv\right] p_{n,k}\left( t\right) dt.
\end{equation*}%
Now we estimate $B_{i,n}g$ for $i=1,2,3,4,$ respectively. Let us estimate $%
B_{1,n}g$ as follows 
\begin{eqnarray*}
\left\vert B_{1,n}g\right\vert &\leq &\frac{1}{2X}\sum\limits_{k=0}^{\left[
nx\right] }\left\vert k-nx\right\vert p_{n,k}\left( x\right) \left(
n+1\right) \int\limits_{0}^{k/n}\left( t-x\right) ^{2}\left\vert
\int\limits_{x}^{t}g^{\prime \prime \prime }\left( v\right) dv\right\vert
p_{n,k}\left( t\right) dt \\
&=&\frac{1}{2X}\sum\limits_{k=0}^{\left[ nx\right] }\left( k-nx\right)
p_{n,k}\left( x\right) \left( n+1\right) \int\limits_{0}^{k/n}\left(
t-x\right) ^{2}\left[ \int\limits_{x}^{t}\left\vert g^{\prime \prime \prime
}\left( v\right) \right\vert dv\right] p_{n,k}\left( t\right) dt \\
&\leq &\frac{1}{2X}\sum\limits_{k=0}^{\left[ nx\right] }\left( k-nx\right)
p_{n,k}\left( x\right) \left( n+1\right) x^{2}\int\limits_{0}^{k/n}\left[
\int\limits_{x}^{t}\left\vert g^{\prime \prime \prime }\left( v\right)
\right\vert dv\right] p_{n,k}\left( t\right) dt \\
&=&\frac{1}{2X}\sum\limits_{k=0}^{\left[ nx\right] }\left( nx-k\right)
p_{n,k}\left( x\right) \left( n+1\right) x^{2}\int\limits_{0}^{k/n}\left[
\int\limits_{t}^{x}\left\vert g^{\prime \prime \prime }\left( v\right)
\right\vert dv\right] p_{n,k}\left( t\right) dt \\
&\leq &\frac{1}{2X}\sum\limits_{k=0}^{\left[ nx\right] }\left( nx-k\right)
p_{n,k}\left( x\right) \left( n+1\right) x^{2}\int\limits_{0}^{k/n}\left[
\int\limits_{0}^{1}\left\vert g^{\prime \prime \prime }\left( v\right)
\right\vert dv\right] p_{n,k}\left( t\right) dt \\
&\leq &\frac{1}{2X}\left\Vert g^{\prime \prime \prime }\right\Vert
\sum\limits_{k=0}^{\left[ nx\right] }\left( nx-k\right) p_{n,k}\left(
x\right) \left( n+1\right) x^{2}\int\limits_{0}^{1}p_{n,k}\left( t\right) dt
\end{eqnarray*}%
\begin{equation*}
=\frac{x}{2\left( 1-x\right) }\left\Vert g^{\prime \prime \prime
}\right\Vert \sum\limits_{k=0}^{\left[ nx\right] }\left( nx-k\right)
p_{n,k}\left( x\right) \leq \frac{x}{2\left( 1-x\right) }\left\Vert
g^{\prime \prime \prime }\right\Vert \left[ nxT_{0,n}\left( x\right)
-T_{1,n}\left( x\right) \right]
\end{equation*}%
where%
\begin{equation*}
\left\Vert f\right\Vert :=\left\Vert f\right\Vert _{L_{1}\left( 0,1\right) }.
\end{equation*}%
In view of (\ref{mom}) we obtain%
\begin{equation}
\left\Vert B_{1,n}g\right\Vert =0.  \label{7}
\end{equation}%
Analogously, $B_{2,n}g$ can be estimated by%
\begin{equation*}
\left\vert B_{2,n}g\right\vert \leq \frac{1}{2X}\sum\limits_{k=0}^{\left[ nx%
\right] }\left\vert k-nx\right\vert p_{n,k}\left( x\right) \left( n+1\right)
\int\limits_{k/n}^{x}\left( t-x\right) ^{2}\left\vert
\int\limits_{x}^{t}\left\vert g^{\prime \prime \prime }\left( v\right)
\right\vert dv\right\vert p_{n,k}\left( t\right) dt
\end{equation*}%
\begin{eqnarray*}
&=&\frac{1}{2X}\sum\limits_{k=0}^{\left[ nx\right] }\left( k-nx\right)
p_{n,k}\left( x\right) \left( n+1\right) \int\limits_{k/n}^{x}\left(
t-x\right) ^{2}\left[ \int\limits_{x}^{t}\left\vert g^{\prime \prime \prime
}\left( v\right) \right\vert dv\right] p_{n,k}\left( t\right) dt \\
&\leq &\frac{1}{2X}\sum\limits_{k=0}^{\left[ nx\right] }\left( k-nx\right)
p_{n,k}\left( x\right) \left( n+1\right) \left( \frac{k}{n}-x\right)
^{2}\int\limits_{k/n}^{x}\left[ \int\limits_{x}^{t}\left\vert g^{\prime
\prime \prime }\left( v\right) \right\vert dv\right] p_{n,k}\left( t\right)
dt \\
&=&\frac{1}{2X}\sum\limits_{k=0}^{\left[ nx\right] }\left( nx-k\right)
p_{n,k}\left( x\right) \left( n+1\right) \left( \frac{k}{n}-x\right)
^{2}\int\limits_{k/n}^{x}\left[ \int\limits_{t}^{x}\left\vert g^{\prime
\prime \prime }\left( v\right) \right\vert dv\right] p_{n,k}\left( t\right)
dt \\
&\leq &\frac{1}{2X}\sum\limits_{k=0}^{\left[ nx\right] }\left( nx-k\right)
p_{n,k}\left( x\right) \left( n+1\right) \left( \frac{k}{n}-x\right)
^{2}\int\limits_{k/n}^{x}\left[ \int\limits_{0}^{1}\left\vert g^{\prime
\prime \prime }\left( v\right) \right\vert dv\right] p_{n,k}\left( t\right)
dt \\
&\leq &\frac{\left\Vert g^{\prime \prime \prime }\right\Vert }{2X}%
\sum\limits_{k=0}^{\left[ nx\right] }\left( nx-k\right) p_{n,k}\left(
x\right) \left( n+1\right) \left( \frac{k}{n}-x\right)
^{2}\int\limits_{k/n}^{x}p_{n,k}\left( t\right) dt \\
&\leq &\frac{\left\Vert g^{\prime \prime \prime }\right\Vert }{2X}%
\sum\limits_{k=0}^{\left[ nx\right] }\left( nx-k\right) p_{n,k}\left(
x\right) \frac{\left( k-nx\right) ^{2}}{n^{2}}=\left\Vert g^{\prime \prime
\prime }\right\Vert \sum\limits_{k=0}^{\left[ nx\right] }\left\vert
k-nx\right\vert p_{n,k}\left( x\right) \frac{\left( k-nx\right) ^{2}}{n^{2}}
\\
&\leq &\frac{\left\Vert g^{\prime \prime \prime }\right\Vert }{2n^{2}X}%
\sum\limits_{k=0}^{n}\left\vert k-nx\right\vert p_{n,k}^{1/2}\left( x\right)
\left( k-nx\right) ^{2}p_{n,k}^{1/2}\left( x\right) .
\end{eqnarray*}%
by using H\"{o}lder inequality and (\ref{mom}) 
\begin{eqnarray*}
&\leq &\frac{\left\Vert g^{\prime \prime \prime }\right\Vert }{2n^{2}X}%
\left( \sum\limits_{k=0}^{n}\left( k-nx\right) ^{2}p_{n,k}\left( x\right)
\right) ^{1/2}\left( \sum\limits_{k=0}^{n}\left( k-nx\right)
^{4}p_{n,k}\left( x\right) \right) ^{1/2} \\
&=&\frac{\left\Vert g^{\prime \prime \prime }\right\Vert }{2n^{2}X}\left(
T_{2,n}\left( x\right) -2nxT_{1,n}\left( x\right) +\left( nx\right)
^{2}T_{0,n}\left( x\right) \right) ^{1/2} \\
&&\times \left( T_{4,n}\left( x\right) -4nxT_{3,n}\left( x\right) +6\left(
nx\right) ^{2}T_{2,n}\left( x\right) -4\left( nx\right) ^{3}T_{1,n}\left(
x\right) +\left( nx\right) ^{4}T_{0,n}\left( x\right) \right) ^{1/2} \\
&=&\frac{\left\Vert g^{\prime \prime \prime }\right\Vert }{2n^{2}X}\left(
nX\right) ^{1/2}\left( 3\left( nX\right) ^{2}+(1-6X)nX\right) ^{1/2} \\
&=&\frac{\left\Vert g^{\prime \prime \prime }\right\Vert }{2n}\left(
3nX+1-6X\right) ^{1/2}
\end{eqnarray*}%
So if we take norm of $B_{2,n}g$, we get 
\begin{equation}
\left\Vert B_{2,n}g\right\Vert \leq \frac{1}{\sqrt{n}}\left\Vert g^{\prime
\prime \prime }\right\Vert .  \label{8}
\end{equation}%
Now we estimate $B_{3,n}g$ with (\ref{mom}) as following%
\begin{equation*}
\left\vert B_{3,n}g\right\vert \leq \frac{1}{2X}\sum\limits_{k=\left[ nx%
\right] +1}^{n}\left\vert k-nx\right\vert p_{n,k}\left( x\right) \left(
n+1\right) \int\limits_{x}^{k/n}\left( t-x\right) ^{2}\left\vert
\int\limits_{x}^{t}\left\vert g^{\prime \prime \prime }\left( v\right)
\right\vert dv\right\vert p_{n,k}\left( t\right) dt
\end{equation*}%
\begin{eqnarray*}
&=&\frac{1}{2X}\sum\limits_{k=\left[ nx\right] +1}^{n}\left( k-nx\right)
p_{n,k}\left( x\right) \left( n+1\right) \left( \frac{k}{n}-x\right)
^{2}\int\limits_{x}^{k/n}\left[ \int\limits_{x}^{t}\left\vert g^{\prime
\prime \prime }\left( v\right) \right\vert dv\right] p_{n,k}\left( t\right)
dt \\
&\leq &\frac{\left\Vert g^{\prime \prime \prime }\right\Vert }{2n^{2}X}%
\sum\limits_{k=0}^{n}\left( k-nx\right) ^{3}p_{n,k}\left( x\right) \\
&=&\frac{\left\Vert g^{\prime \prime \prime }\right\Vert }{2n^{2}X}\left[
T_{3,n}\left( x\right) -3nxT_{2,n}\left( x\right) +3\left( nx\right)
^{2}T_{1,n}\left( x\right) -\left( nx\right) ^{3}T_{0,n}\left( x\right) %
\right] \\
&=&\frac{\left\Vert g^{\prime \prime \prime }\right\Vert }{2n}(1-2X).
\end{eqnarray*}%
So, we get the following inequality for the norm of $B_{3,n}g$%
\begin{equation}
\left\Vert B_{3,n}g\right\Vert \leq \frac{1}{2n}\left\Vert g^{\prime \prime
\prime }\right\Vert .  \label{9}
\end{equation}%
Finally, we estimate $B_{4,n}g$ with the (\ref{mom}) 
\begin{equation*}
\left\vert B_{4,n}g\right\vert \leq \frac{1}{2X}\sum\limits_{k=\left[ nx%
\right] +1}^{n}\left\vert k-nx\right\vert p_{n,k}\left( x\right) \left(
n+1\right) \int\limits_{k/n}^{1}\left( t-x\right) ^{2}\left\vert
\int\limits_{x}^{t}\left\vert g^{\prime \prime \prime }\left( v\right)
\right\vert dv\right\vert p_{n,k}\left( t\right) dt
\end{equation*}%
\begin{eqnarray*}
&=&\frac{1}{2X}\sum\limits_{k=\left[ nx\right] +1}^{n}\left( k-nx\right)
p_{n,k}\left( x\right) \left( n+1\right) \left( 1-x\right)
^{2}\int\limits_{k/n}^{1}\left[ \int\limits_{0}^{1}\left\vert g^{\prime
\prime \prime }\left( v\right) \right\vert dv\right] p_{n,k}\left( t\right)
dt \\
&\leq &\frac{1}{2X}\left\Vert g^{\prime \prime \prime }\right\Vert
\sum\limits_{k=\left[ nx\right] +1}^{n}\left( k-nx\right) p_{n,k}\left(
x\right) \left( 1-x\right) ^{2} \\
&\leq &\frac{1-x}{2x}\left\Vert g^{\prime \prime \prime }\right\Vert \left[
T_{1,n}\left( x\right) -nxT_{0,n}\left( x\right) \right] =0,
\end{eqnarray*}%
and hence%
\begin{equation}
\left\Vert B_{4,n}g\right\Vert =0.  \label{10}
\end{equation}%
Collecting the results in (\ref{7})-(\ref{10}), we have by (\ref{12})%
\begin{equation*}
\left\Vert R_{n}g\right\Vert \leq \frac{1}{\sqrt{n}}\left\Vert g^{\prime
\prime \prime }\right\Vert +\frac{1}{2n}\left\Vert g^{\prime \prime \prime
}\right\Vert \leq \frac{2}{\sqrt{n}}\left\Vert g^{\prime \prime \prime
}\right\Vert .
\end{equation*}%
Finally we obtain by using (\ref{11})%
\begin{equation*}
\left\Vert \left( D_{n}g\right) ^{\prime }-g^{\prime }\right\Vert \leq \frac{%
2}{n+2}\left\Vert g^{\prime }\right\Vert +\frac{2}{n+2}\left\Vert g^{\prime
\prime }\right\Vert +\frac{2}{\sqrt{n}}\left\Vert g^{\prime \prime \prime
}\right\Vert .
\end{equation*}%
According to Stein's inequality (see, e.g., \cite{13}, Theorem A10.1) one has%
\begin{eqnarray*}
\left\Vert g^{\prime \prime }\left( x\right) \right\Vert _{L_{1}\left(
0,1\right) } &\leq &C\sqrt{\left\Vert g^{\prime }\left( x\right) \right\Vert
_{L_{1}\left( 0,1\right) }\left\Vert g^{\prime \prime \prime }\right\Vert
_{L_{1}\left( 0,1\right) }} \\
&\leq &C\left( \left\Vert g^{\prime }\right\Vert _{L_{1}\left( 0,1\right)
}+\left\Vert g^{\prime \prime \prime }\right\Vert _{L_{1}\left( 0,1\right)
}\right) ,
\end{eqnarray*}%
where $C>1$. For $\frac{1}{n+2}<\frac{1}{\sqrt{n}}$ we have%
\begin{equation*}
\left\Vert \left( D_{n}g\right) ^{\prime }-g^{\prime }\right\Vert \leq \frac{%
2(C+1)}{\sqrt{n}}\left( \left\Vert g^{\prime }\right\Vert +\left\Vert
g^{\prime \prime \prime }\right\Vert \right) .
\end{equation*}%
So the proof is now complete.


\begin{thebibliography}{99}
\bibitem{1} J. L. Durrmeyer, Une formule d'inversion de la transform%
%TCIMACRO{\U{b4}}%
%BeginExpansion
\'{}%
%EndExpansion
ee de Laplace: Applications`a la th%
%TCIMACRO{\U{b4}}%
%BeginExpansion
\'{}%
%EndExpansion
eorie des moments, Th`ese de 3e Cycle, Facult%
%TCIMACRO{\U{b4}}%
%BeginExpansion
\'{}%
%EndExpansion
e des Sciences del'Universit%
%TCIMACRO{\U{b4}}%
%BeginExpansion
\'{}%
%EndExpansion
e de Paris, 1967.

\bibitem{2} A. Lupas, Die Folge der Betaoperatoren, Dissertation,
Universitat Stuttgart, (1972)

\bibitem{3} M.M. Derriennic, Sur l'approximation de fonctions int6grables
sur [0, I] par des polyn6mes de Bernsteinmodifies, J. Approx. Theory 31,
325-343, (1981).

\bibitem{4} H.H. Gonska and X. Zhou, A global inverse theorem on
simultaneous approximation by Bernstein-Durrmeyeroperators, d. Approx.
Theory 67, 284-302, (1991).

\bibitem{5} Z. Ditzian and K. Ivanov, Bernstein-type operators and their
derivatives, J. Approx. Theory 56, 72-90,(1989).

\bibitem{6} Bardaro, C.,Butzer, P.L., Stens, R. L., and Vinti, G.,
Convergence in Variation and Rates of Approximation for Bernstein-Type
Polynomials and Singular Convolution \.{I}ntegrals, Analysis (Munich), 23
(4), 299-346 (2003).

\bibitem{7} Lorentz, G. G., Bernstein polynomials, University of Toronto
Press, Toronto (1953).

\bibitem{8} Agratini O., On the variation detracting property of a class of
operators, Apply. Math. Lett. 19, 1261-1264, (2006).

\bibitem{9} P. Pych-Taberska and H. Karsli, On the rates of convergence of
Bernstein-Chlodovsky polynomials and their B\'{e}zier-type variants.
Applicable Analysis 90(3-4),403-416 (2011).

\bibitem{10} Kivinukk, A., Metsmagi, T.: Approximation in variation by the
Meyer--K\"{o}nig and Zeller operators. Proc. Estonian Acad. Sci.60(2),
88--97 (2011)

\bibitem{11} \"{O}. \"{O}ks\"{u}zer, H. Karsli, F. Tasdelen, On Convergence
of Bernstein-Stancu Polynomials\ in the Variation Seminorm, Numer. Funct.
Anal. Optim. 37(4), 1-20 (2016).

\bibitem{12} H. G\"{u}l \.{I}nce \.{I}larslan, G. Ba\c{s}canbaz-Tunca,
Convergence in Variation for Bernstein-Type Operators, Mediterr. J. Math.
DOI 10.1007/s00009-015-0640-1 (2015).

\bibitem{13} Trigub, R. M. and Belinsky, E. S. Fourier Analysis and
Approximation of Functions. Kluwer Academic Publishers, Dordrecht (2004).
\end{thebibliography}
\end{document}